\input amstex
\input amsppt.sty

\define\su{\sum\limits}
\define\pro{\prod\limits}
\define\inte{\int\limits}

\NoRunningHeads
\TagsOnRight
\magnification 1200
\topmatter
\title
Point Processes and the\\ Infinite Symmetric Group \\Part IV:
Matrix Whittaker kernel
\endtitle
\author
Alexei Borodin
\endauthor
\abstract
We study a 2--parametric family of probability measures on the space
of countable point configurations on the punctured real line (the points of the random configuration are concentrated near zero). These
measures (or, equivalently, point processes) have been introduced in
Part II (A.~Borodin, math/9804087) in connection with the problem
of harmonic analysis on the infinite symmetric group. The main result of the present paper is a determinantal formula for the correlation
functions. 

The formula involves a kernel called the matrix Whittaker kernel.
Each of its two diagonal blocks governs the projection of the
process on one of the two half--lines; the corresponding kernel on the half--line was studied in Part III (A.~Borodin and G.~Olshanski, math/9804088). 

While the diagonal blocks of the matrix Whitaker kernel are symmetric, the whole kernel turns
out to be $J$-symmetric, i.e., symmetric with respect to a natural indefinite inner product.  

We also discuss a rather surprising connection of our processes with the recent work by B.~Eynard and M.~L.~Mehta (cond-mat/9710230) on correlations of eigenvalues of coupled random matrices.
\endabstract
\thanks Supported by
the Russian Program for Support of Scientific Schools under grant
96-15-96060. 
\endthanks
\endtopmatter
\document
\head Introduction\endhead
In this paper we continue the study of a remarkable family of stochastic
point processes started in [O1] (Part I), [B1] (Part II), and [BO1] (Part III).
This Part IV is followed by the next paper [O2] (Part V) which concludes our
work on the subject. The summary of results from Parts I--V is presented in
[BO2].

The point processes in question live on the punctured interval
$I=[-1,1]\setminus\{0\}$. They are parametrized by two parameters $z$ and
$z'$ which satisfy one of the following two conditions:
$$
\align
(*) \quad \ \,&z'=\bar{z},\ z\in\Bbb C\setminus\Bbb Z;              \\
(**)\quad &z\text{ and } z'\text{ are real and for a certain }m\in\Bbb Z,\
m< z,z'< m+1.
\endalign
$$
 We denote the corresponding process by $\Cal
P_{zz'}$. The definition of these processes can be found in [O1, \S4] and
[B1, Introduction].

The origin of our study is the work [KOV] on generalized regular
representations of the infinite symmetric group --- the processes $\Cal
P_{zz'}$ govern, in a certain sense, the decomposition of these
representations into irreducibles.

Part I was devoted to the translation of the representation theoretic problem
in the language of the point processes, and to the calculation of the density
function of the processes. In Part II we computed the higher correlation
functions of $\Cal P_{zz'}$ and introduced a modification of the processes
(so--called `lifting') which substantially simplified the formulas for the
correlation functions. Part III contains several applications of the results
of Parts I,II.

In Part II we gave explicit integral representations for all higher
correlation functions, see [B1, Theorem 1.7.3, Theorem 2.2.1]. However, these
integral representations were rather complicated. After lifting, see [B1,
Chapter 3], multidimensional integrals involved in the integral
representations were reduced to two--dimensional ones. Moreover, we have
proved that in the domains where all variables have the same sign, the
correlation functions are given by determinantal formulas with a certain
kernel on $\Bbb R_+$ called the Whittaker kernel [B, Theorem 3.3.4].

In this paper we show that after lifting the correlation functions are
everywhere given by the determinantal formulas with a certain kernel
defined on $\Bbb R^*$. It is convenient to write this kernel in the matrix
form according to the splitting $\Bbb R^*=\Bbb R_+\sqcup\Bbb R_-$. We call
the
kernel thus obtained the {\it matrix Whittaker kernel}. One of its blocks
coincides with the Whittaker kernel introduced in [B1].

The expression of the `lifted' correlation functions via the matrix Whittaker
kernel is our main result.

Determinantal form for the correlation functions appears in many problems of
random matrix theory and mathematical physics, see, e.g., [Dy], [Me2], [KBI].
In most situations the kernels are symmetric or hermitian (see, however,
[B2]).
But the matrix Whittaker kernel turns out to be `$J$-symmetric' where
$J=\left[\matrix 1&0\\0&-1\endmatrix\right]$, see Remark 2.9 below.
The appearance of kernels with such symmetry in determinantal formulas
for correlation functions seems to be new.

At the end of the paper we consider the systems of eigenvalues of two random
coupled matrices studied in [IZ], [Me1], [MS], [EM], [MN], [Ey]. As was
recently proved in [EM], the correlation functions of such systems are also
given by determinantal formulas. We show that this result and our
considerations have common combinatorial background. Here we give an
alternative proof of the main result of [EM] in case of two coupled matrices
([EM] deals with a more general situation of several matrices coupled in a
chain). Our proof seems to be more straightforward and less elegant than
that from [EM].

The paper is organized as fol
lows. In Section 1 we prove that the correlation
functions of $\Cal P_{zz'}$ after lifting are given by determinantal formulas
and express the kernel through two--dimensional distributions introduced
in [B1]. The proof is purely combinatorial. In Section 2 we show that the
kernel obtained in Section 1 can be conveniently written via the Whittaker
functions. Section 3 deals with the eigenvalue correlations of two random
coupled matrices.

I am very greatful to G.~I.~Olshanski for his permanent attention and
support. His `jewel' interpretation of the proof of Proposition 1.1 made the
presentation much more transparent.

 \head 1. Determinantal formula \endhead
Let us introduce the notation used in [B1].

Denote by $\Phi_{n,d}$ the set of mappings
$$
\varphi :\{1,\ldots,n\}\to \{1,1';\ldots;d,d'\}
$$
subject to the two conditions

1) $\varphi$ is injective, i.e. $\varphi(i)\ne \varphi(j)$ if $i\ne j$;

2) $\operatorname{Im} \varphi \cap\{m,m'\}\ne\emptyset$ for all
$m=1,\ldots,d.$

It is clear that $\Phi_{n,d} \ne \emptyset$ if and only if $n/2\le d\le n$.

Starting from a function (or a distribution) $F(r_1,s_1;\ldots;r_d,s_d)$ in $2d$
 variables we define the function $(\varphi F)(x_1,\ldots, x_n)$ in $n$ variables,
  $\varphi \in \Phi_{n,d}$, as follows.
Let us rename the variable $r_i$ by $x_k$ if $\varphi (k)=i$, $s_j$ by $x_k$ if
$\varphi (k)=j'$, and let us do this for all $r_i,s_j$ such that $i,j'\in
 \operatorname{Im} \varphi$. Then let us integrate $F$ over all $r_l,s_m$ such
  that
$l \notin \operatorname{Im} \varphi, \ m'\notin \operatorname{Im}
\varphi$. The result is denoted by $(\varphi F)(x_1,\ldots, x_n)$.

Let $N(x,y)$ and $w(x,y)$ be (generalized) functions in two variables such that
$$\operatorname{supp}N(x,y)\subset\{(x,y)\ |x\geq 0,y\leq 0\},
$$
$$
\operatorname{supp}w(x,y)\subset\{(x,y)\ |x\leq 0,y\geq 0\}.
$$
For all $d=1,2,\dots$ set
$$
H_d(r_1,s_1;\ldots;r_d,s_d)=\frac
1{d!}\prod_{i=1}^dw(s_i,r_i)\cdot\det[N(r_i,s_j)]_{i,j=1}^d.
\tag 1.1
$$
Our main result in this section is the following statement.
\proclaim{Proposition 1.1}
Let $x_1,\ldots,x_n\neq 0$. Then
$$
\sum_{d}\sum_{\varphi\in\Phi_{n,d}}(\varphi
H_d)(x_1,\dots,x_n)=\det[\widetilde{N}(x_i,x_j)]_{i,j=1}^n
\tag 1.2
$$
where
$$
\widetilde{N}(x,y)=\cases
\int\limits_sN(x,s)w(s,y)ds,&x,y>0\\
N(x,y),&x>0,\ y<0\\
\iint\limits_{r,s}w(x,r)N(r,s)w(s,y)drds-w(x,y),&x<0,\ y>0\\
\int\limits_r w(x,r)N(r,y) dr,&x,y<0
\endcases
$$
\endproclaim
\example{Remark 1.2} Theorem 3.2.6 of [B1]
 claims that the $n$th lifted correlation function of our processes is equal to the
left-hand side of (1.2) for
the following choice of $N(x,y)$ and $w(x,y)$:
$$
N(x,-y)=t\left(\phi_{z'}(x)\phi_{-z'}(y)e^{-x-y}\right)\odot \left(\phi_{z,-z-
1}(x)\phi_{-z,z-1}(y)\right)\odot\left(\delta(x-y)\chi_{[0,1]}(x)\right),
\tag 1.3
$$
$$
w(x,y)=\frac{1}{|x|+y},
\tag 1.4
$$
see Section 1.4 and 1.6 of [B1] for the notation.

(Note that our function $N$ here differs from $N$ defined by [B1, (3.6)] by
changing the sign of $y$ and multiplication by the constant $t=zz'$.)
\endexample
\demo{Proof of Proposition 1.1} The proof will consist of five steps.

{\bf Step 1.} On the set $\Phi_{n,d}$ there is a natural action of the
symmetric group $S_d$ which permutes the pairs $(i,i')$ for $i=1,\ldots,d$.
On the other hand, the function $H_d$ defined in (1.1) is invariant with
respect to the permutations of the pairs $(r_i,s_i)$, $i=1,\ldots,d$. This
implies that all elements $\varphi\in \Phi_{n,d}$ lying in the same orbit of
$S_d$ give the same contribution to the left-hand side of (1.2). Each orbit
contains $d!$ elements, hence, we may assume that the summation in the
left-hand side of (1.2) is taken over the set of orbits of $S_d$ in
$\Phi_{n,d}$ and throw out the factor $\frac 1{d!}$ in the right-hand side of
(1.1), cf. [B1, Proposition 1.7.1].

{\bf Step 2.} Let us plug the formula (1.1) into the left-hand side of (1.2) and
expand all the determinants. Then the summands will be parametrized by the
triples $(d,\widehat\varphi,\sigma)$ where $\widehat\varphi$ is an orbit of
$S_d$ in $\Phi_{n,d}$ and $\sigma\in S_d$ parametrizes the expansion of the
determinant (1.1). We may also consider $\sigma$ as a map from
$\{1,2,\ldots,d\}$ to the set $\{1',2',\ldots,d'\}$.

We associate to each such summand an oriented (not necessarily connected)
necklace with $2d$ black and white beads, some of the beads are marked by the
numbers from $1$ to $n$, in the following way.

Let us fix a representative $\varphi\in \widehat\varphi\subset \Phi_{n,d}$.

The necklace contains $d$ white and $d$ black beads, every two neighbors are
of different colors. White beads correspond to the elements $\{1,2,\ldots,d\}$
and black ones -- to the elements $\{1',2',\ldots,d'\}$ of the set
$\{1,1';\ldots;d,d'\}$ used in the definition of $\Phi_{n,d}$.

A white bead $i$ is marked by the number $k$ if $\varphi(k)=i$; a black bead
$j'$ is marked by the number $l$ if $\varphi(l)=j'$.

The white bead $i$ always follows the black bead $i'$ and is followed by the
black bead $\sigma(i)$.

Now we forget that the beads correspond to the elements of $\{1,1';\ldots;d,d'\}$
and obtain an oriented necklace with black and white beads some of which are
marked by $1,2,\ldots,n$.

It is not difficult to see that different representatives of the same orbit
$\widehat\varphi$ produce the same necklace, and different orbits correspond
to different necklaces.

For a fixed necklace the corresponding summand of the left-hand side of (1.2)
can be constructed as follows. We take $2d$ indeterminates and associate them
to the beads somehow. Then for each edge (interval between two beads) we
produce a factor like this:

 for an edge going from a white bead (indeterminate $r$) to a black bead
(indeterminate $s$) we take $N(r,s)$;

 for an edge going from a black bead (indeterminate $s$) to a black bead
(indeterminate $r$) we take $w(s,r)$.

 Thus, we obtain $d$ factors, we multiply them and attribute to the
indeterminates corresponding to the beads marked by $1,\ldots,n$ the names
$x_1,\ldots,x_n$ respectively.

 Finally, we integrate over all other (mute) indeterminates and add the sign
$\operatorname{sgn} \sigma$. The sign will be discussed in Step 5.

{\bf Step 3.} Now we shall translate to the `jewel language' the right-hand side
of (1.2).

We may assume that all variables $x_1,\ldots,x_n$ have definite signs.

First, we expand the determinant in the right-hand side of (1.2), the summands
are parametrized by permutations $\tau\in S_n$. Further, if a summand
contains a factor $\widetilde{N}(x,y)$ with negative $x$ and positive $y$, then
$$
\widetilde{N}(x,y)= \iint\limits_{r,s}w(x,r)N(r,s)w(s,y)drds-w(x,y),
\tag 1.5
$$
and we break this summand into two: one will contain
$$
\iint\limits_{r,s}w(x,r)N(r,s)w(s,y)drds,
\tag 1.6
$$
and the other
$$
-w(x,y).\tag 1.7
$$

For each summand we produce an oriented (not necessarily connected) necklace
with $n$ beads corresponding to the variables $x_1,\ldots,x_n$; the beads
corresponding to positive variables are white, to negative -- black. All beads are
marked by the numbers $1,\ldots,n$ according to the numbers of corresponding
variables. The bead marked by $i$ is followed by the bead $\tau(i)$.

On each edge going from a black (`negative') bead $i$ to a white (`positive')
bead $j$ we put the sign `$N$' or `$w$' according to which of the summands (1.6)
or (1.7) was chosen while breaking $N(x_i,x_j)$ into two parts.

In the first case we shall say that we used $N$-choice, in the second -- $w$-
choice.

{\bf Step 4.} Now we shall establish one-to-one corresponds
between the sets of necklaces introduced in Steps 2 and 3.

First, let us take a necklace from Step 2. It
contains $2d$
black and white beads some of which are marked by the numbers
$1,\ldots,n$. Then we can immediately produce an oriented
necklace with $n$ beads by throwing out all mute (unmarked)
beads.

The second condition in the definition of $\Phi_{n,d}$ implies
that at least one end of each edge going from black to white is not mute (at least
one element in each pair $(i,i')$ is covered by $\varphi\in \Phi_{n,d}$). It means
that after removing all mute beads, we can obtain an edge `black--to--white' in
two different ways: it either comes from the initial necklace unchanged, or there
were two mute beads removed from the middle of this edge. In the first case we
put the sign `$w$' on this edge, in the second -- the sign `$N$'.

Thus, starting with a necklace from Step 2 we produce a necklace from Step 3.
Let us go the other way around.

We take a necklace with $n$ beads from Step 3. Then we follow the rules:

(1) inside any edge of the type `white--to--white' we insert one mute black bead,

(2) inside any edge of the type `black--to--black' we insert one mute white bead,

(3) the edges `white--to--black' remain unchanged,

(4) the edges `black--to--white' with the sign `$w$' remain unchanged,

(5) inside any edge  `black--to--white' with the sign `$N$' we insert one white
and one black mute beads, white goes first.

Clearly, two procedures described above establish one--to--one correspondence
between the necklaces of Steps 2 and 3, thus, between the summands of the left--
hand side and the right--hand side of (1.2).

Straightforward check shows that up to sign the corresponding summands from
Steps 2 and 3 are equal. To compute, up to sign, the summand corresponding to
a pair of necklaces, we take each edge of the necklace with $n$ beads and
according to 5 different situations for this edge described above (suppose the
edge goes from the $i$th bead to the $j$th bead) we produce the factor:
$$
\align
&(1):\int_sN(x_i,s)w(s,x_j)ds;\\
&(2):\int_r w(x_i,r)N(r,x_j) dr;\\
&(3):N(x_i,x_j);\\
&(4):w(x_j,x_i);\\
&(5):\iint_{r,s}w(x_i,r)N(r,s)w(s,x_j)drds-w(x_i,x_j).
\endalign
$$

Then all such $n$ factors are multiplied.

{\bf Step 5.} Let us verify that corresponding summands have the same sign.

We take two corresponding necklaces: one from Step 2, the other from Step 3.

The sign for the first one equals $\operatorname{sgn} \sigma$ ($\sigma\in S_d$
was introduced in Step 2). The number of cycles of $\sigma$ is equal to the
number of connected components of the necklace, let us denote it by $c$. Then
$\operatorname{sgn}\sigma=(-1)^{d-c}.$

The sign for the second necklace is equal to $\operatorname{sgn} \tau$ ($\tau$
was introduced in Step 3) times $-1$ to the power of the number of $w$--
choices, because of the sign `$-$' in (1.7).
The number of cycles in $\tau$ is equal to $c$, because both necklaces have the
same number of connected components. Thus,
$\operatorname{sgn}\tau=(-1)^{n-c}.$
Furthermore, the number of $w$--choices is equal to $n-d$, it can be easily
verified by induction on the number of $w$--choices.
Hence
$$
n-c-\text{number of $w$--choices}=d-c,
$$
and the signs are equal. \qed
\enddemo
\head 2. The matrix Whittaker kernel\endhead
Our goal in this section is to compute explicitly the kernel $\widetilde N(x,y)$
introduced in Proposition 1.1 for the special $N(x,y)$ and $w(x,y)$ given by
(1.3) and (1.4).

It is convenient to write $\widetilde N(x,y)$ in the matrix form
$$
\widetilde N(x,y) =\left[\matrix
N_{++}(x,y)&N_{+-}(x,y)\\N_{-+}(x,y)&N_{--}(x,y)\endmatrix\right]
$$
where
$$
\alignat{3}
N_{++}(x,y)&=\int\limits_sN(x,s)w(s,y)ds, &x,y>0;\qquad\\
N_{+-}(x,-y)&=N(x,y),&x>0,y<0;\qquad\\
N_{-+}(-x,y)&=\iint\limits_{r,s}w(x,r)N(r,s)w(s,y)drds-
w(x,y),\qquad&x<0,y>0;\qquad\\
N_{--}(-x,-y)&=\int\limits_t  w(x,r)N(r,y)dr,&x,y<0.\qquad
\endalignat
$$
Starting from now we set, see Remark 1.2,
$$
N(x,-y)=t\left(\phi_{z'}(x)\phi_{-z'}(y)e^{-x-y}\right)\odot \left(\phi_{z,-z-
1}(x)\phi_{-z,z-1}(y)\right)\odot\left(\delta(x-y)\chi_{[0,1]}(x)\right),
\tag 2.1
$$
$$
w(x,y)=\cases\frac{1}{|x|+y},& x<0,y>0\\
0,&\text{otherwise}
\endcases.
\tag 2.2
$$

Moreover, in all our calculations below we shall assume that
$$
-1<\Re z,\,\Re z'<0.
\tag 2.3
$$
Under this assumption several distributions used below become just integrable
functions, which will be suitable for us. Our final result will depend on $z$ and
$z'$ analytically. Thus, by analytic continuation, it will hold for all pairs $z,z'$.

We use the standard notation $W_{\kappa,\mu}(x)$ for  {\it the Whittaker
functions}, see [E, Chapter 6].
Note that
$$
W_{\kappa,\mu}(x)=W_{\kappa,-\mu}(x).
$$

\proclaim{Proposition 2.1}
$$
\gathered
N_{+-}(x,y)=N(x,-y)=\frac{\sin\pi z\sin\pi z'}{\pi^2}\left(\frac
xy\right)^{\frac{z+z'}2}
e^{\frac{-x-y}2}(xy)^{-\frac 12}\\ \times \frac{W_{\frac{z+z'+1}{2},\pm\frac{z-
z'}2}(x)W_{\frac{-z-z'+1}2,
  \pm\frac{z-z'}2}(y)
   +tW_{\frac{z+z'-1}2,\pm\frac{z-z'}2}(x)W_{\frac{-z-z'-1}2,
  \pm\frac{z-z'}2}(y)}{x+y}.
\endgathered
\tag 2.4
$$
\endproclaim
\demo{Proof} Let us multiply (2.4) by $(x+y)$ and compare the
$(k,l)$-moment of both sides. We shall need the relations
$$
\langle\phi_a(x)e^{-x},x^m\rangle=(a+1)_m,\quad
  \langle\phi_{a,b}(x),x^m\rangle=\frac{(a+1)_m}{\Gamma(a+b+m+2)}
$$
which are just Euler gamma and beta integrals; and
$$
\int\limits_0^\infty
x^{b-\frac c2-1}e^{-\frac x2}W_{\frac c2-a,\frac c2 -\frac
12}(x)dx=\frac{\Gamma(b)\Gamma(b-c+1)}{\Gamma(a+b-c+1)},
$$
which is a version of [E, 6.10(7)].

The moments of a pseudoconvolution of distributions are the products of
moments
of the factors. Thus, the $(k,l)$-moment of $(x+y)N(x,y)$ equals
$$
\align
\langle(x+y)N(x,-y),x^ky^l\rangle&=t\,\frac{(z+1)_{k+1}(z'+1)_{k+1}(-z+1)_l(-
z'+1)_l}
{(k+1)!l!(k+l+2)}\\&+t\,\frac{(z+1)_{k}(z'+1)_{k}(-z+1)_{l+1}(-z'+1)_{l+1}}
{k!(l+1)!(k+l+2)}\\&=
t\,\frac{(z+1)_{k}(z'+1)_{k}(-z+1)_l(-z'+1)_l}
{k!l!}\\&+t^2\frac{(z+1)_{k}(z'+1)_{k}(-z+1)_l(-z'+1)_l}
{(k+1)!(l+1)!}.
\endalign
$$

Straightforward computation shows that these two summands are exactly the
$(k,l)$-moments of two summands in the numerator of the right-hand side of
(2.4).
Thus, we proved that the moments of both sides of (2.4) coincide after
multiplication by $(x+y)$.

In general, it is not sufficient to conclude that our distributions coincide,
because they are not compactly supported. However,
it is easy to check that both sides of (2.4) after multiplication by $(x+y)$
can be realized as liftings of some compactly
supported distributions.(See Chapter 3 of [B1] for the definition of the
lifting.) By comparison of moments, these compactly supported distributions
must coincide, and so do their liftings. \qed \enddemo \proclaim{Proposition
2.2} $$ \gathered
N_{++}(x,y)=\int\limits_s\frac{N(x,-
s)}{s+y}ds=\frac{1}{\Gamma(z)\Gamma(z')}
\left(\frac xy\right)^{\frac{z+z'}2}e^{\frac {y-x}2}(xy)^{-\frac 12}\\ \times
\frac{W_{\frac{z+z'+1}{2},\pm\frac{z-z'}2}(x)W_{\frac{z+z'-1}2,
  \pm\frac{z-z'}2}(y)
   -W_{\frac{z+z'-1}2,\pm\frac{z-z'}2}(x)W_{\frac{z+z'+1}2,
  \pm\frac{z-z'}2}(y)}{x-y}
\endgathered
\tag 2.5
$$
\endproclaim
\demo{Proof}
 This claim is implicitly proved in [B1]. Let us explain why.

Application of the Stieltjes transform to a pseudoconvolution is equivalent to
the
application of it to one of the factors, if the others does not have a
singularity at zero. Assuming that $-1<\Re z,\Re z'<0$, see (2.3) above, we may
take the Stieltjes
transform of the second factor of (2.1), instead of applying it to the whole
$N(x,y)$.

 Then we use Lemma 2.2.3 from [B1] which claims that
$$
\inte_s \frac{\phi_{-z,z-1}(s)\,ds}{s+y}=y^{-z}(1+y)^{z-1}.
\tag 2.6
$$
We get
$$
\multline
N_{++}(x,y)=\int\limits_s\frac{N(x,-s)}{s+y}ds=t\left(\phi_{z'}(x)\phi_{-
z'}(y)e^{-x-y}\right) \\\odot \left(\phi_{z,-z-1}(x)\,y^{-z}(1+y)^{z-
1}\right)\odot\left(\delta(x-y)\chi_{[0,1]}(x)\right).
\endmultline
$$
By definition of pseudoconvolution (see [B1, 1.4]), we get
$$
\aligned
N_{++}(x,y)=
t\int\limits_{u,v}
\phi_{z'}(u)&\phi_{-z'}(v)e^{-u-v}\\ \times
\int\limits_{ 0}^1 &\phi_{z,-z-1}\left(\frac{x}{uw}\right)
\left(\frac{y}{vw}\right)^{-z}
\Bigl(1+\frac{y}{vw}\Bigr)^{z-1}
\frac{dwdudv}{w^2uv},
\endaligned
\tag 2.7
$$
cf. [B1, (2.6)].

Lemma 2.2.4 of [B1] states that
$$
\inte_0^1 \phi_{z,-z-1}\left(\frac{r_1}{uw}\right)\left(\frac{r_2}{vw}\right)^{-z}
\left(1+\frac{r_2}{vw}\right)^{z-1}
\frac{dw}{uvw^2}=\frac{\phi_{-z}(u-r_1)\phi_z(r_2+v)}{r_1v+r_2u}.
$$

Applying this to (2.7) we get
$$
N_{++}(x,y)=t\int\limits_{u,v}
\phi_{z'}(u)\phi_{-z'}(v)
\phi_{-z}(u-x)\phi_{z}(v+y)\frac{e^{-u-v}}{xu+yv}dudv.
$$
By introducing new integration
variables
$$
t_1=u/x-1,\quad t_2=v/y
$$
we see that
$$
\multline
N_{++}(x,y)=t\,\left(\frac xy\right)^{z'}e^{\frac{y-x}2}\\
\times\iint\limits_{t_1,t_2}\phi_{-z}(t_1)\phi_{-
z'}(t_2)\phi_{z'}(t_1+1)\phi_{z}(t_2+1)\frac{e^{-x(t_1+1/2)-
y(t_2+1/2)}}{t_1+t_2+1}dt_1dt_2,
\endmultline
$$
cf. [B1, Theorem 3.3.1].
The formulas (3.8) and (3.9) from [B1] (together with the definition of the
kernel $M(x,y)$, see [B1, Theorem 3.3.1]) conclude the proof.\qed \enddemo By
changing the signs of the parameters $z$ and $z'$ in Proposition 2.2 and
changing $x\leftrightarrow y$, we get the following statement.
\proclaim{Corollary 2.3}
$$
\gathered
N_{--}(x,y)=\int\limits_r\frac{N(r,-y)}{r+x}dr=\frac{1}{\Gamma(-z)\Gamma(-
z')}
\left(\frac xy\right)^{\frac{z+z'}2}e^{\frac{x-y}2}(xy)^{-\frac 12}\\ \times
\frac{W_{\frac{-z-z'+1}2,\pm\frac{z-z'}2}(x)W_{\frac{-z-z'-1}2,
  \pm\frac{z-z'}2}(y)-W_{\frac{-z-z'-1}{2},\pm\frac{z-z'}2}(x)W_{\frac{-z-z'+1}2,
  \pm\frac{z-z'}2}(y)}{x-y}.
\endgathered
\tag 2.8
$$
\endproclaim
Finally, we compute $N_{-+}(x,y)$.
\proclaim{Proposition 2.4}
$$
\gathered
N_{-+}(x,y)=
\iint\limits_{r,s}\frac{N(r,-s)}{(r+x)(s+y)}drds-
\frac 1{x+y}=-\left(\frac xy\right)^{\frac{z+z'}2}e^{\frac{x+y}2}
(xy)^{-\frac 12}\\ \times \frac{W_{\frac{-z-z'+1}{2},\pm\frac{z-z'}2}(x)
W_{\frac{z+z'+1}2,
  \pm\frac{z-z'}2}(y)
   +tW_{\frac{-z-z'-1}2,\pm\frac{z-z'}2}(x)W_{\frac{z+z'-1}2,
  \pm\frac{z-z'}2}(y)}{x+y}.
\endgathered
\tag 2.9
$$
\endproclaim
\example{Remark 2.5}
Note that
$$
N_{-+}(x,y)=-\frac{\pi^2}{\sin \pi z\sin\pi z'}\, e^{x+y}N_{+-}(y,x).
\tag 2.10
$$
This is a non--trivial and surprising fact. A detailed discussion of this
`coincidence' (or rather of its modification, see Remark 2.9) can be found in
[O2].
\endexample
\demo{Proof of Proposition 2.4} First, as in the
 proof of Proposition 2.2, we apply both Stieltjes transforms to the second
factor of the pseudoconvolution (2.1). Using (2.6) twice, we get $$ \multline
\int\limits_{r,s}\frac{N(r,-s)}{(r+x)(s+y)}drds=t\left(\phi_{z'}(x)\phi_{-z'}(y)e^{-
x-y}\right) \\\odot \left(x^{z}(1+x)^{-z-1}\,y^{-z}(1+y)^{z-
1}\right)\odot\left(\delta(x-y)\chi_{[0,1]}(x)\right).
\endmultline
$$
By definition of the pseudoconvolution, we get
$$
\multline\int\limits_{r,s}\frac{N(r,-s)}{(r+x)(s+y)}drds=
t\int\limits_{u,v}
\phi_{z'}(u)\phi_{-z'}(v)e^{-u-v}\\ \times
\int\limits_{ 0}^1
\left(\frac{x}{uw}\right)^{z}
\Bigl(1+\frac{x}{uw}\Bigr)^{-z-1}
\left(\frac{y}{vw}\right)^{-z}
\Bigl(1+\frac{y}{vw}\Bigr)^{z-1}
\frac{dwdudv}{w^2uv}.
\endmultline
$$
Introducing new variables
$$
\cases \tau_1=u/x&\\
 \tau_2= v/y&
\endcases
$$
we obtain
$$
\multline
\int\limits_{r,s}\frac{N(r,-s)}{(r+x)(s+y)}drds=
t\left(\frac xy\right)^{z'}\int\limits_{\tau_1,\tau_2}
\phi_{z'}(\tau_1)\phi_{-z'}(\tau_2)e^{-x\tau_1-y\tau_2}\\ \times
\int\limits_{ 0}^1
(1+\tau_1w)^{-z-1}(1+\tau_2w)^{z-1}dwd\tau_1d\tau_2.
\endmultline
$$
We shall need the following lemma.
\proclaim{Lemma 2.6}
$$
\int\limits_{ 0}^1
(1+\tau_1w)^{-z-1}(1+\tau_2w)^{z-1}dw=\frac
1z\frac1{\tau_2-\tau_1}\left(\left(\frac{1+\tau_2}{1+\tau_1}\right)^z-
1\right)
$$
\endproclaim
\demo{Proof of Lemma 2.6} First we use [E, 5.8.2(5)]:
$$
\int\limits_{ 0}^1
(1+\tau_1w)^{-z-1}(1+\tau_2w)^{z-1}dw=F_1(1,z+1,-z+1,2;-\tau_1-\tau_2),
$$
where $F_1(\alpha,\beta,\beta',\gamma;x,y)$ is the Appell's hypergeometric
function. By [E, 5.10(1)], we can simplify the last expression:
$$
F_1(1,z+1,-z+1,2;-\tau_1-\tau_2)=(1+\tau_2)^{-1}
F\left(1,z+1;2;\frac{\tau_2-\tau_1}{\tau_2+1}\right).
$$

But
$$
F(1,z+1;2;\tau)=\frac 1{z\tau}((1-\tau)^{-z}-1),
$$
and we arrive at our claim.
\qed
\enddemo
Applying the lemma,
$$
\aligned\int\limits_{r,s}\frac{N(r,-s)}{(r+x)(s+y)}drds=
t&\left(\frac xy\right)^{z'}\int\limits_{\tau_1,\tau_2}
\phi_{z'}(\tau_1)\phi_{-z'}(\tau_2)e^{-x\tau_1-y\tau_2}\\ &\times
\frac
1z\,\frac1{\tau_2-\tau_1}\left(\left(\frac{1+\tau_2}{1+\tau_1}\right)^z-
1\right)
d\tau_1d\tau_2.
\endaligned
\tag 2.11
$$
Let us multiply this relation by $x+y$. Note that
$$
(x+y)e^{-x\tau_1-y\tau_2}=-\left(\frac\partial{\partial\tau_1}+
\frac\partial{\partial\tau_2}\right)e^{-x\tau_1-y\tau_2},
$$
cf. proof of [B1, Theorem 3.3.4].

 We are going to integrate (2.11) multiplied by $x+y$ by
parts. The following formula is easy to check
$$
\multline
\left(\frac\partial{\partial\tau_1}+
\frac\partial{\partial\tau_2}\right)\left[\phi_{z'}(\tau_1)\phi_{-z'}(\tau_2)
\,\frac
1z\,\frac1{\tau_2-\tau_1}\left(\left(\frac{1+\tau_2}{1+\tau_1}\right)^z-
1\right)\right]\\ =\frac 1{t}\,\phi_{z'-1}(\tau_1)\phi_{-z'-1}(\tau_2)
-\frac 1{t}\,\phi_{z'-1}(\tau_1)\phi_{-z'-1}(\tau_2)(1+\tau_1)^{-z}({1+\tau_2})^z
\\-
\phi_{z'}(\tau_1)\phi_{-z'}(\tau_2)(1+\tau_1)^{-z-1}(1+\tau_2)^{-z-1}.
\endmultline
$$
Hence, integrating by parts,
$$
\multline
(x+y)\int\limits_{r,s}\frac{N(r,-s)}{(r+x)(s+y)}drds=
\left(\frac xy\right)^{z'}\int\limits_{\tau_1,\tau_2} \biggl(
\phi_{z'-1}(\tau_1)\phi_{-z'-1}(\tau_2)
\\-\phi_{z'-1}(\tau_1)\phi_{-z'-1}(\tau_2)({1+\tau_1})^{-z}
({1+\tau_2})^{z}\\-
t\,\phi_{z'}(\tau_1)\phi_{-z'}(\tau_2)(1+\tau_1)^{-z-1}(1+\tau_2)^{z-1}
\biggr)e^{-x\tau_1-y\tau_2}d\tau_1d\tau_2.
\endmultline
$$
Using Euler gamma integral for the first summand and the standard integral
representation of the
Whittaker functions
$$
W_{\kappa,\mu}(x)=e^{-x/2}x^{\mu+1/2}\int\limits_{\tau}
\phi_{\mu-\kappa-1/2}(\tau)(1+\tau)^{\mu+\kappa-1/2}e^{-\tau x}d\tau
$$
for the last two,
we rewrite the last relation in the form
$$
\multline
(x+y)\int\limits_{r,s}\frac{N(r,-s)}{(r+x)(s+y)}drds=
1-\left(\frac xy\right)^{\frac{z+z'}2}e^{\frac{x+y}2}(xy)^{-\frac12}\\
 \times\biggl(W_{\frac{-z-z'+1}{2},\frac{z'-z}2}(x)
W_{\frac{z+z'+1}2,
  \frac{z-z'}2}(y)
   +tW_{\frac{-z-z'-1}2,\frac{z'-z}2}(x)W_{\frac{z+z'-1}2,
  \frac{z-z'}2}(y)\biggr)
\endmultline
$$
which is equivalent to (2.9). \qed
\enddemo
Note that in all four quadrants the expression for $\widetilde N(x,y)$ has the
same factor
$$
\left|\frac xy\right|^{\frac{z+z'}2}e^{\frac {y-x}2}=
\frac {|x|^{\frac{z+z'}{2}}e^{-\frac x2}}{|y|^{\frac{z+z'}{2}}e^{-\frac y2}}.
\tag 2.12
$$
(For example, if $x,y<0$, $\widetilde N(-x,-y)=N_{--}(x,y)$ contains the factor
$$
\left(\frac xy\right)^{\frac{z+z'}2}e^{\frac{x-y}2}
$$
which after changing the signs of $x$ and $y$ coincides (2.12).)

But all factors (2.12) disappear in the determinants of the form $\det\widetilde
N(x_i,x_j)$. Thus, we can introduce a new kernel as follows.

Denote
$$
\align
A_+(x)&=x^{-\frac12}W_{\frac{z+z'+1}{2},\pm\frac{z-z'}2}(x),\\
B_+(x)&= x^{-\frac12}W_{\frac{z+z'-1}{2},\pm\frac{z-z'}2}(x),\\
A_-(x)&=x^{-\frac12}W_{\frac{-z-z'+1}{2},\pm\frac{z-z'}2}(x),\\
B_-(x)&=x^{-\frac12}W_{\frac{-z-z'-1}{2},\pm\frac{z-z'}2}(x).
\endalign
$$
For $x,y>0$ set
$$
\alignat{2}
K(x,y)=K_{++}(x,y)&=\frac{1}{\Gamma(z)\Gamma(z')}\,\frac{A_+(x)
B_+(y)-B_+(x)A_+(y)}{x-y}\\
K(x,-y)=K_{+-}(x,y)&=\frac{\sqrt{\sin\pi z\sin\pi z'}}{\pi}\,
\frac{A_+(x) A_-(y)+t\,B_+(x) B_-(y)}{x+y}\\
K(-x,y)=K_{-+}(x,y)&=-\frac{\sqrt{\sin\pi z\sin\pi z'}}{\pi}\,\frac{A_+(y)
A_-(x)+t\, B_+(y)B_-(x)}{x+y}\\
K(-x,-y)=K_{--}(x,&y)=\frac{1}{\Gamma(-z)\Gamma(-z')}\,\frac{A_-(x)
B_-(y)-B_-(x)A_-(y)}{x-y}
\endalignat
$$
or, in matrix form,
$$
K(x,y) =\left[\matrix
K_{++}(x,y)&K_{+-}(x,y)\\K_{-+}(x,y)&K_{--}(x,y)\endmatrix\right].
$$
We call $K(x,y)$ the {\it matrix Whittaker kernel}.

All our previous work is summarized by the following statement.
\proclaim{Theorem 2.7} The $n$th correlation function of the lifted process
$\widetilde {\Cal P_{zz'}}$ has the form
$$
\tilde\rho_n^{(zz')}(x_1,\ldots,x_n)=\det[K(x_i,x_j)]_{i,j=1}^n
$$
where $K(x,y)$ is the matrix Whittaker kernel.
\endproclaim
\example{Remark 2.8} Theorem 3.3.4 of [B1] where we established the
determinantal form of the lifted correlation functions on a part of the phase
spase is a direct corollary of Theorem 2.7.  The Whittaker kernel introduced
in [B1] is the `++' block of the matrix Whittaker kernel.
\endexample
\example{Remark 2.9} Note that all blocks of the matrix Whittaker kernel are
real--valued functions and
 $$
\gather
 K_{++}(x,y)=K_{++}(y,x),\quad
K_{--}(x,y)=K_{--}(y,x),\\
 K_{+-}(x,y)=-K_{-+}(y,x),
\endgather
 $$
cf.  Remark 2.5.
This means that the matrix Whittaker kernel is not symmetric but
`$J$--symmetric' for $$J=\left[\matrix 1&0\\0&-1\endmatrix\right].  $$ A
general discussion of stochastic point processes governed by such kernels can
be found in [O2].  \endexample

 \head 3.  Two random coupled matrices
\endhead In this section we will explain the connection of our problems with
the system of eigenvalues of two random coupled matrices.  This system was
introduced in [IZ] and studied in [Me1], [MS], [EM], [MN], [Ey]. In
particular, B.~Eynard and M.~L.~Mehta proved in [EM] that the correlation
functions of this system are given by determinantal formulas (in fact, they
proved the determinantal formulas for the correlation functions in case of
finitely many matrices coupled in a chain).  Let us describe this
result in more details.  We will follow [EM] slightly changing the notation.

Consider two  complex hermitian $N \times N$ matrices $A$ and
$B$ with the
probability density
$$
F(A,B)=const\cdot \exp(-\operatorname{tr}\{U(A)+V(B)+cAB\})
$$
where $U(x)$ and $V(x)$ are real polynomials of even degree
with positive coefficients of their highest powers and $c$ is a
real constant.

Let us denote the sets of (real) eigenvalues of $A$ and $B$ by
$\{x_{11},\ldots,x_{1N}\}$ and $\{x_{21},\ldots,x_{2N}\}$. As was
proved in [IZ], [Me1], the probability density for the eigenvalues has
the form
$$
\aligned
p_N(x_{11},\ldots,x_{1N};&x_{21},\ldots,x_{2N})\\=const\cdot &\det
[w(x_{2i},x_{1j})]_{i,j=1}^N\cdot\prod_{i<j}[(x_{1i}-x_{1j})(x_{2i}-x_{2j})]
\endaligned
\tag 3.1
$$
where
$$
w(x,y)=\exp(-U(y)-V(x)+cxy).
$$
The correlation functions are defined as follows:
$$
\aligned
\rho_{k,l}(x_{11},\ldots,x_{1k};x_{21},\ldots,x_{2l})&\\=\frac{N!^2}{(N-
k)!\,(N-l)!}\int
p_N&(x_1;x_2)\,dx_{1k+1}\ldots dx_{1N}\,dx_{2l+1}\ldots dx_{2N}.
\endaligned
\tag 3.2
$$
Under certain non-degeneracy conditions on the pairing
$$
\langle f(x),g(y)\rangle=\inte_{x,y}f(x)g(y)w(y,x)dxdy,
\tag 3.3
$$
 we can find two systems of biorthogonal polynomials $$\{P_0(x),\ldots,P_{N-
1}(x)\}\  \text{and} \ \{Q_0(y),\ldots,Q_{N-1}(y)\}$$ with respect to this pairing.
In other words, $\deg P_i=\deg Q_i=i$ and
$$
\langle P_i,Q_j\rangle=\delta_{ij}.
$$

Set
$$
H(x,y)=\su_{i=0}^{N-1}P_i(x)Q_i(y)
\tag 3.4
$$

\proclaim{Proposition 3.1 ([EM])}
$$
\rho_{k,l}(x_{11},\ldots,x_{1k};x_{21},\ldots,x_{2l})=\det
[K_{ij}(x_{ir},x_{js})]_{i,j=1,2;r=1,\ldots,k;s=1,\ldots,l}
\tag 3.5
$$
where
$$
\aligned
K_{11}(x,y)&=\int\limits_sH(x,s)w(s,y)ds;\\
K_{12}(x,y)&=H(x,y);\\
K_{21}(x,y)&=\iint\limits_{r,s}w(x,r)H(r,s)w(s,y)drds-w(x,y);\\
K_{22}(x,y)&=\int\limits_t w(x,r)H(r,y)dr.
\endaligned
\tag 3.6
$$
\endproclaim

\example{Remark 3.2}
Note that the formula (3.6) is identical to the formula for
the kernel $\widetilde N(x,y)$ in Proposition 1.1. This is not an accidental
coincidence. It turns out that there exists a wide class of measures on the
infinite--dimensional Thoma simplex such that after lifting the correlation
functions of these measures are given by determinantal formulas with a matrix
kernel. Moreover, all these matrix kernels have the form (3.6) for an
appropriate distribution $H(x,y)$ (e.g., the matrix Whittaker kernel).

In particular, this class of measures include many finite--dimensional measures
for which the techniques of [EM] can be applied. For example, one of such
measures after lifting has the density of the form (3.1) with the weight function
 $w(x,y)=\frac 1{x+y}$ on $\Bbb R_+\times \Bbb R_+$, cf. (2.2).
\endexample

Now we shall give a sketch of the proof of Proposition 3.1 which is different
from the proof given in [EM]. Namely, we shall reduce the statement to
Proposition 1.1 proved in Section 1.
\demo{Sketch of the proof of Proposition 3.1}
Following the notation of Section 1, let us introduce the set  $\Phi_{k,l;d}$ of
mappings
$$
\varphi:\{{}_11,\ldots,{}_1k;{}_21,\ldots,{}_2l\}\to \{1,1';\ldots,d,d'\}
$$
subject to the three conditions

1) $\varphi$ is injective;

2) $\operatorname{Im} \varphi \cap\{m,m'\}\ne\emptyset$ for all
$m=1,\ldots,d$;

3) $\varphi({}_{1}m)\in\{1,\ldots,d\},$ $\varphi({}_{2}m)\in\{1',\ldots,d'\}$ for
all $m$.

By analogy with Section 1, for each $\varphi\in\Phi_{k,l;d}$ and a function
$F(r_1,s_1;\ldots;r_d,s_d)$ in $2d$
variables we define the function $(\varphi
F)(x_{11},\ldots,x_{1k};x_{21},\ldots,x_{2l})$ in $k+l$ variables. Specifically,
according to $\varphi$ we attribute to some variables $r_i$ names $x_{1m}$, to
some variable $s_j$ the names $x_{2m}$, and then integrate $F$ over all mute
variables.

We shall prove that
$$
\rho_{k,l}(x_{11},\ldots,x_{1k};x_{21},\ldots,x_{2l})=\sum_{d}\sum_{\varphi\in\Phi_{k,l;d}}(\varphi H_d)(x_{11},\ldots,x_{1k};x_{21},\ldots,x_{2l})
\tag 3.7
$$
for
$$
H_d(r_1,s_1;\ldots;r_d,s_d)=\frac
1{d!}\prod_{i=1}^dw(s_i,r_i)\cdot\det[H(r_i,s_j)]_{i,j=1}^d.
\tag 3.8
$$
Then word--for--word repetition of the `jewel' proof of Proposition 1.1 will
prove (3.5); the only difference is that in this case we consider the
necklaces where white beads are marked by the elements of the set
$\{{}_11,\ldots,{}_1k\}$ and black beads are marked by the elements of the
set $\{{}_21,\ldots,{}_2l\}$.

In order to prove (3.7), we shall first show that the sets $\Phi_{k,l;d}$
naturally appear in the expansion of the determinant of a matrix along $k$
rows and $l$ columns and then apply this observation to the
determinant in the RHS of (3.1).

Let $M=(m_{ij})_{i,j=1}^N$ be a sufficiently large square matrix. For any
$\varphi \in \Phi_{k,l;d}$ we define $(\varphi M)$ as a certain sum of
expressions of the form
$$
\pm m_{a_1b_1}\cdots m_{a_db_d}\,
M\binom{\cdots\hat a_1\cdots\hat a_d\cdots}{\cdots\hat b_1\cdots \hat
b_d\cdots} \tag 3.9 $$ where by $M\binom{\cdots\hat a_1\cdots\hat
a_d\cdots}{\cdots\hat b_1\cdots \hat b_d\cdots}$ we denote the determinant of
the submatrix of $M$ obtained from $M$ by removing the rows
$a_1,\ldots,a_d$
and the columns $b_1,\ldots,b_d$.

The sum is defined as follows. We associate the indices
 $\{a_1,\ldots,a_d;
b_1,\ldots,b_d\}$
 with the elements of
 $\{1,\ldots,d;1',\dots,d'\}$
and the numbers $\{1,\ldots,k;1,\ldots,l\}$ with the elements of
$\{{}_11,\ldots,{}_1k;{}_21,\ldots,{}_2l\}$.
 Then, according to our map
$$
\varphi:
\{{}_11,\ldots,{}_1k;{}_21,\ldots,{}_2l\}
\to
 \{1,\ldots,d;1',\dots,d'\}
$$
we specify some of $a_i$'s and $b_i$'s. As for the rest, we let $d-k$ other
(not specified) $a_i$'s vary over all $(d-k)$--tuples of pairwise distinct
numbers from the set $\{k+1,\ldots,N\}$, and $d-l$ other $b_i$'s vary over
all $(d-l)$--tuples of pairwise distinct numbers from $\{l+1,\ldots,N\}$.

The sign in (3.9) is chosen in such a way that (3.9) enters the expansion of
$\det M$ with positive sign.

\proclaim{Lemma 3.3}
For any integer $k,l\geq 0$ and sufficiently large square matrix $M$
$$
\det M=\su_d\su_{\varphi\in\Phi_{k,l;d}}\frac 1{d!}(\varphi M).
\tag 3.10
$$
\endproclaim
This formula is the expansion of $\det M$ along the first $k$ rows and $l$
columns. For example, for $k=1$, $l=0$, (3.10) turns into the usual first row
expansion of $\det M$:
$$
\det M=\su_{i=1}^N (-1)^{i+1}m_{1i}\cdot M\binom{\hat 1\cdots}{\cdots \hat
i\cdots}.  $$
For $k=1,\ l=1$ we get the formula
$$ \det M=m_{11}\cdot M\binom {\hat
1\cdots}{\hat 1\cdots}+\su_{i,j=2}^N (-1)^{i+j+1}m_{i1}m_{1j}\cdot M\binom{\hat
1\cdots\hat i\cdots}{\hat 1\cdots\hat j\cdots}.  $$ Here the first term
 corresponds to the unique element of $\Phi_{1,1;1}$ and the second term
corresponds to the two elements of $\Phi_{1,1;2}$, which give the same
contribution.

 For arbitrary $k$ and $l=0$, (3.10) coincides
with the well--known Laplace expansion of the determinant.

The proof of Lemma 3.3 is straightforward. The factor $\frac 1{d!}$ in (3.10)
may be removed if the summation is taken over the sets of orbits of $S_d$'s
in $\Phi_{k,l;d}$'s, cf. Step 1 of the proof of Proposition 1.1.

Using elementary row and column transformations in the Vandermonde
determinants of (3.1), we can rewrite (3.1) in the form
$$
\aligned p_N(x_{11},\ldots,x_{1N};&\
x_{21},\ldots,x_{2N})\\=\frac{1}{N!^2}\cdot &
\det
[w(x_{2i},x_{1j})]_{i,j=1}^N\cdot\det
[P_i(x_{1j})]_{i,j=1}^N\cdot[Q_i(x_{2j})]_{i,j=1}^N.
\endaligned \tag 3.11
$$

Now we apply (3.10) to the first determinant of (3.11). We shall show that
after normalization and integration over extra variables as in (3.2), the
$\varphi$--term of the expansion of (3.11) will give exactly the
$\varphi$--term of the RHS of (3.7) for all $\varphi$'s.

So, we set $M=(w(x_{2i},x_{1j}))_{i,j=1}^N$ and employ Lemma 3.3.

Note that for a fixed $\varphi$ the contributions of all summands (3.9) to
the RHS of (3.2) are equal because they differ by permutations of the
integration variables $x_{1k+1},\dots,x_{1N}$ and $x_{2l+1},\dots,x_{2N}$.
The total number of summands (3.9) is the number of $(d-k)$--tuples and
$(d-l)$--tuples of pairwise distinct numbers taken from the sets with $N-k$
and $N-l$ elements respectively. This number equals
$$
\frac{(N-k)!}{(N-d)!}\,\frac{(N-l)!}{(N-d)!}.
\tag 3.12
$$
So we keep this combinatorial factor and choose one suitable for us summand
of the form (3.9), namely, such that
$$
\gather
\{a_1,\ldots,a_d\}=\{1,\dots,d\}, \\
\{b_1,\ldots,b_d\}=\{1,\dots,d\}.
\endgather
$$
Then this summand has the form
$$
\frac{\operatorname{sgn}\tau}{N!^2}\pro_{i=1}^dw(x_{2\tau(i)},x_{1i})\cdot
\det[w(x_{2i},x_{1j})]_{i,j=d+1}^N\cdot\det
[P_i(x_{1j})]_{i,j=1}^N\cdot[Q_i(x_{2j})]_{i,j=1}^N.
\tag 3.13
$$
for a certain $\tau\in S_d$.

The expression (3.13) can be easily integrated over $x_{1d+1},\dots x_{1N};
x_{2d+1},\dots,x_{2N}$. By expanding all the determinants, using orthogonality
conditions and Gram's formula
we get
$$
\frac{(N-d)!^2}{N!^2}\,\pro_{i=1}^dw(x_{2\tau(i)},x_{1i})\cdot
\det[H(x_{1i},x_{2\tau(i)})]_{i,j=1}^d
\tag 3.14
$$
where $H(x,y)$ is defined in (3.4). Finally, integrating (3.14) over
$x_{1k+1},\dots,x_{1d}$ and $x_{2l+1},\ldots x_{2d}$ and multiplying it by
(3.12), combinatorial factor from (3.2) and $\frac 1{d!}$ from (3.10), we
arrive at $(\varphi H_d)(x_{11},\dots,x_{1k};x_{21},\ldots,x_{2l})$.\qed
\enddemo

\bigskip
\bigskip

{\smc A.~Borodin}: Department of Mathematics, The University of
Pennsylvania, Philadelphia, PA 19104-6395, U.S.A.  E-mail address:
{\tt borodine\@math.upenn.edu}

\enddocument